\documentclass[12pt]{article}
\usepackage{amsxtra,amssymb,amsthm,amsmath,latexsym}

\theoremstyle{plain}
\def\R{{\mathbb R}}

\def\C{{\mathbb C}}

\def\oH{{\overset{\circ}{H}}}
\def\oH1{{\overset{\circ}{H}\kern-.02in{}^1}}

\def\bee{\begin{equation*}}
\def\eee{\end{equation*}}
\def\be{\begin{equation}}
\def\ee{\end{equation}}

\begin{document}

%\begin{titlepage}
\title{When is $F(p)$  the Laplace transform of a bounded $f(t)$?
}

\author{Alexander G. Ramm\\
 Department  of Mathematics, Kansas State University, \\
 Manhattan, KS 66506, USA\\
ramm@ksu.edu\\
%,\\ %fax 785-532-0546, tel. 785-532-0580}
http://www.math.ksu.edu/\,$\sim$\,ramm}
%}

\date{}
\maketitle\thispagestyle{empty}

%%%%%%%%%%%%%%%%%%%%%%%%%%%%%%%%%%%%%%%%%%%%%%%%%
\begin{abstract}

Sufficient conditions are given for a function $F(p)$, analytic in Re$p>0$,
to be a Laplace transform   of a function $f(t)$, such that 
$max_{t\ge 0}|f(t)|<\infty$, $f(0)=0$.  

\end{abstract}
%%%%%%%%%%%%%%%%%%%%%%%%%%%%%%%%%%%%%%%%%%%%%%%%%
%\end{titlepage}

\section{Introduction }\label{S:1}
 There is a large literature on the Laplace transform, e.g., \cite{S}, \cite{LS}.
 Tables of the Laplace transform of distributions are published, see, e.g., \cite{BP}. The Laplace transform of a function $f(t)$, $f(t)=0$ for $t<0$,
  $|f(t)|\le Ce^{at}$ for some positive constants $C$ and $a$,
  is defined by the formula
\be\label{e1}
L(f):=F(p):=\int_0^\infty e^{-pt}f(t)dt, \quad p=s+i\eta,\,\, s> a.
\ee
   Under this assumption 
 $F(p)$ is an analytic function of $p$ in the half-plane $ s>a$ and $|F(p)|\le \frac  C{s-a}$,
$s=Re\, p$, as follows from the estimate 
\be\label{e1b}
|F(p)|=|\int_0^\infty e^{-pt}f(t)dt|\le C\int_0^\infty e^{-(s-a)t}dt=\frac C {s-a}, \quad s>a.
\ee
 We assume that
 \be\label{e1a}
 \lim_{|p|\to \infty, -\frac {\pi} 2\le \phi\le \frac {\pi} 2} F(p)=0,
 \ee
 where $\phi=arg\, p$ is the argument of $p$.
  
  If   $f\in L^2(0,\infty)$, then $a=0$. In Section 3 we discuss possible generalizations.
   
  We are interested in the sufficient conditions for
   $F(p)$ to be the Laplace transform of a function $f(t)$, such that
\be\label{e2}
\sup_{t\ge 0}|f(t)|<\infty, \quad f(0)=0.
\ee
These conditions are \eqref{e1b}, \eqref{e1a} and
\be\label{e3}
|F(i\eta)|=O(|i\eta|^{-b}),\,\,b>1\text {  for  } |\eta|\gg 1; \,\, F(p) \text {is analytic for }  s>0.  
\ee
The condition $F(i\eta)\in L^2(-\infty, \infty)$ follows from the assumption $f(t)\in L^2(0, \infty)$  by the Parceval's identity. This
condition is not necessary for $F(p)$ to be the Laplace transform of $f(t)$ satisfying conditions \eqref{e2}.

Our result is formulated in the following theorem.
 
 {\bf Theorem 1.} {\em If conditions  \eqref{e1b}, \eqref{e1a} and \eqref{e3} hold, then $F(p)=L(f)$ and
 $f(t)$ satisfies conditions  \eqref{e2}.}
  
  This result is new. It differs from the known results:
  usually some assumptions are made on $f(t)$ and {\em necessary} conditions  are derived for $F(p)$ to be the Laplace transform. In Theorem 1 assumptions are made on $F(p)$ and  {\em sufficient} conditions are given  for $F(p)$ to be the Laplace transform  of a function 
  $f(t)$ satisfying conditions \eqref{e2} and vanishing for $t<0$.
    
 \section{Proof}\label{S:2}

Let
\be\label{e4}
f (t) =\frac 1 {2\pi} \int_{-\infty}^{\infty} e^{it\eta} F(i\eta)d\eta.
\ee

   The integral \eqref{e4}  converges since  $b>1$. Therefore, $L(f ) = F(p)$, $f (t)$ is a continuous function and $\sup_{t\ge 0}|f (t)|<
\infty$. This follows from formula \eqref{e4} and the known inversion formula for the Laplace transform,  \cite{LS}, namely:
\be\label{e4a}
f(t)=\frac 1 {2\pi i}\int_{C_\sigma}e^{pt}F(p)dp,
\ee
where $C_\sigma$, $\sigma>a$, is the straight line $\sigma-i\eta, \sigma+i\eta$, $-\infty<\eta<\infty$. In our case $\sigma=0$,  $dp=id\eta$
and $C_\sigma$ is the straight line $-i\eta, i\eta$,  $-\infty<\eta<\infty$.

For convenience of the reader, let us give a version of a proof of the inversion formula. Under our assumptions $f(t)$, defined in \eqref{e4} is a continuous
function on $[0,\infty)$, uniformly bounded because the function $F(i\eta)$  
is absolutely integrable on  the whole axis $-\infty<\eta< \infty$. 
Consider 
\be\label{e5a} 
\int_0^\infty f(t)e^{-qt}dt=\frac 1 {2\pi} \int_{-\infty}^{\infty} d\eta
F(i\eta)\int_0^\infty e^{-(q-i\eta)t}dt=-\frac 1 {2\pi i} \int_{-i\infty}^{i\infty}F(p)\frac 1 {p-q}dp,\,\,\, p=i\eta. 
\ee
Let $C_n$ be the closed contour consisting of $[in, -in]\cup L_n$, where
$L_n$ is the semi-circle $p=ne^{i\phi}$, $-\frac {\pi}2\le \phi\le \frac {\pi}2$.
Since $b>1$, it follows that 
\be\label{e5b} 
\lim_{n\to \infty}\int_{L_n}F(p)\frac 1 {p-q}dp=0.
\ee
Therefore, the integral on the right side of \eqref{e5a} can be considered
as the integral over the closed contour $C_n=$  with $n\to \infty$. The minus sign in \eqref{e5a} is used to get the closed contour passed counterclockwise.
Consequently, the integral in the right side of \eqref{e5a} is equal to $F(q)$ by the Cauchy formula and the analyticity of $F(p)$ in the half-plane Re$p>0$. Thus,
\be\label{e5c} 
\int_0^\infty f(t)e^{-qt}dt=F(q).
\ee
 Therefore, the inversion formula for the Laplace transform is proved. 
 
Let us prove that $f(0)=0$.
It follows from \eqref{e4} that 
\be\label{e5} 
 f(0)=\frac 1 {2\pi}\int_{-\infty}^{\infty} F(i\eta)d\eta.
 \ee
 If $N$ is sufficiently large for the estimate 
\be\label{e6} 
 |F(p)|\le c(1+|p|)^{-b}, \quad |p|>N,\quad b>1,
 \ee
 to hold, then
 \be\label{e7} 
 f(0)=\frac 1{2\pi}\int_{-N}^{N} F(i\eta)d\eta +o(1) \text{  as  } N\to \infty.
 \ee
 The function $f(t)$ is continuous and uniformly bounded on $[0,\infty)$. Indeed, by the inversion formula \eqref{e4a} with $\sigma=0$, where the integral over the straight line $\sigma=0$ absolutely converges if $b>1$, one has
\be\label{e7a} 
 f(t)=\frac 1{2\pi i}\int_{-iN}^{iN}e^{pt}F(p)dp +\frac 1{2\pi i}\int_{-i\infty}^{-iN}e^{pt}F(p)dp +\frac 1{2\pi i}\int_{iN}^{i\infty}e^{pt}F(p)dp.
 \ee
 The first integral is a continuous function of $t$ because it is taken
 over a compact set and $F(p)\in L^1(-iN,iN)$, the second and third integrals
 are continuous functions of $t$ because $b>1$. At $t=0$ formula \eqref{e7a} reduces to \eqref{e7}.
 
 If we prove that
  \be\label{e8} 
 I_N:=\frac 1{2\pi}\int_{-N}^{N} F(i\eta)d\eta \to 0 \text{  as  } N\to \infty,
 \ee
 then relations \eqref{e7}-\eqref{e8} imply that $f(0)=0$ and Theorem 1 is proved.
 
 To prove \eqref{e8}, consider a closed contour $C_N$, consisting of $[-iN, iN]$ and the semi-circle $L_N:=Ne^{i\phi}$, where $-\frac {\pi}2\le \phi\le  \frac {\pi}2$. Since $F(p)$ is analytic in the half-plane $s>0$,  one has: 
 \be\label{e9} 
 \int_{C_N}F(p)dp=0.
 \ee
  Since 
  $\lim_{N\to \infty}\frac 1{2\pi}\int_{L_N}F(p)dp=0$ because $b>1$, it 
 follows from \eqref{e9} that  equation \eqref{e8} is valid. 
 
 Finally, the condition $f(t)=0$ for $t<0$ follows from the 
  following argument. If $t<0$ then the function $e^{-pt} F(p),\,\,t>0$ is analytic in Re$p>0$  and   
 \be\label{e9a}
 f(-t)=\lim_{N\to \infty}\frac 1{2\pi i}\int_{C_N}e^{-pt}F(p)dp=0, \quad t>0,
 \ee
 because by the Cauchy theorem the integral over the closed contour $C_N$,
 inside which the function $e^{-pt}F(p)$ is analytic, is equal to zero. 
 
  Theorem 1 is proved. \hfill$\Box$
  
The Jordan lemma in the following form, compare \cite{LS}, pp. 412, 469, is useful: 
 
 {\bf Lemma.} {\em If a function $h(p)$ tends to zero uniformly with respect to the argument $\phi$ of $p$, $-\frac {\pi} 2\le \phi\le\frac {\pi} 2$ on the contour $L_n$ as $n\to \infty$ then for $t>0$ one has $\lim_{n\to \infty}\int_{L_n}h(p)e^{pt}dp=0$.}
 % If $t<0$, then $\lim_{n\to \infty}\int_{C^*_n}h(p)e^{pt}dp=0$, where $C^*_n$ is the mirror image of $C_n$ with respect to imaginary axis of the complex $p$-plane.} 

 \section{Discussion }\label{S:3}
In this Section we discuss possible generalization and applications of Theorem 1. We may replace the assumption $F(i\eta)\in L^2(-\infty, \infty)$
by the assumption
\be\label{e10}  
F(i\eta)\in L^p(-iN,iN),\,\,p\ge 1,\quad \text {  \eqref{e3} holds.  }
\ee
The proof of Theorem 1 remains unchanged because the assumption $F(i\eta)\in L^p(-iN,iN),\,\,p\ge 1$  implies that the first integral in \eqref{e7a} is a continuous function of $t$ and assumption \eqref{e6} guarantees that the second and third integrals in \eqref{e7a} are continuous functions of $t$, which
tend to zero as $N\to \infty$. 

In applications Theorem 1 is useful in the definition of hyper-singular integrals, see \cite{R707}. To give an idea of this application, 
consider the integral equation:
\be\label{e11}
q(t)+\frac 1 {\Gamma (\lambda)}\int_0^t (t-\tau)^{\lambda -1}q(\tau)d\tau=f(t).
\ee
If Re$\lambda>0$, then the integral in \eqref{e11} is defined classically,
i.e., from the classical point of view. If $\lambda<0$, then this integral is
a hyper-singular integral,{\em  it diverges classically}. For Re$\lambda>0$ we take the Laplace transform of \eqref{e11} and get
\be\label{e12}
L(q)(1+p^{-\lambda})=L(f),\quad L(q)=\frac{L(f)}{1+p^{-\lambda}},
\ee
where the following formula was used: $L(t^{\lambda-1})=\Gamma (\lambda) p^{-\lambda}$. This formula is valid for all $\lambda\in \C$ except for $\lambda= 0,-1,-2,...$.  The $L(q)$ in \eqref{e12} admits analytic continuation with respect to $\lambda$ from the region Re$\lambda>0$ to the region Re$\lambda<0$, for example, to the point $\lambda=-\frac 1 4$ which is of interest in the Navier-Stokes problem. Theorem 1 is of use to prove that $\frac{L(f)}{1+p^{-\lambda}}$  for $\lambda=-\frac 1 4$ is the Laplace transform of a function $q(t)$ satisfying conditions \eqref{e2}. This can be checked if
 $f(t)$ smooth and rapidly decaying as $t\to \infty$,
so that $|L(f)|\le c(1+|p|)^{-1}$. In this case the function $(1+p^{\frac 1 4})^{-1}$ is analytic in the half-plane Re$p\ge 0$ and is $O(\frac 1 {|p|^{1/4}})$ for $|p|\gg 1$, the function $ L(f)$ is analytic in the half-plane Re$p\ge 0$ and is $O(\frac 1 {|p|})$ for $|p|\gg 1$ on the imaginary axis of the complex plane $p=s+i\eta$. By Theorem 1, the function
$\frac{L(f)}{1+p^{\frac 1 4}}$
is the Laplace transform of a function $q(t)$ satisfying  \eqref{e2}. In this example $b=\frac 5 4>1$.
We have proved the following result.

{\bf Theorem 2.} {\em Assume $\lambda=-\frac 1 4$ and $f(t)$ be a smooth rapidly decaying as $t\to \infty$. Then equation \eqref{e11}
has a unique solution $q(t)$ satisfying \eqref{e2}.
}

The kernel of equation \eqref{e11} with $\lambda=-\frac 1 4$ is
hyper-singular. The integral in this equation with $\lambda=-\frac 1 4$ {\em diverges} classically. Theorem 2 is of prime interest in a study of the Navier-Stokes problem, see \cite{R707}.

 \section{Conclusion }\label{S:4}
 
 Sufficient conditions are given for a function $F(p)$, $p=s+i\eta$, analytic in the half-plane $s>0$  to be the
 Laplace transform of a function $f(t)$, $f(t)=0$ for $t<0$, $\sup_{t\ge 0}|f(t)<\infty$ and $f(0)=0$. This result is useful in a study of the Navier-Stokes problem in $\R^3$ and in a study of integral equations with hyper-singular kernels, see \cite{R707}.

\end{document}